\documentclass[12pt]{amsart}
\usepackage{amssymb}
\title[Special subvarieties of $\sA_g$]
{Special subvarieties of $\sA_g$}
\author[Eckart Viehweg]{Eckart Viehweg}
\address{Universit\"at Duisburg-Essen, Mathematik, 45117 Essen, Germany}
\email{viehweg@uni-essen.de}
\thanks{Both authors are supported by the ``DFG-Schwerpunktprogramm
Globale Methoden in der Komplexen Geometrie'', and the first named author by the DFG-Leibniz program.\\
This work has been achieved thanks to the support
of the European Commission through its 6th Framework Programme
"Structuring the European Research Area" and the contract Nr.
RITA-CT-2004-505493 for the provision of Transnational Access
implemented as Specific Support Action.}
\author[Kang Zuo]{Kang Zuo}
\address{Universit\"at Mainz,
Fachbereich 17, Mathematik,
55099 Mainz, Germany}
\email{kzuo@mathematik.uni-mainz.de}
\begin{document}
\theoremstyle{plain}
\newtheorem{thm}{Theorem}[section]
\newtheorem{theorem}[thm]{Theorem}
\newtheorem{lemma}[thm]{Lemma}
\newtheorem{corollary}[thm]{Corollary}
\newtheorem{proposition}[thm]{Proposition}
\newtheorem{addendum}[thm]{Addendum}
\newtheorem{variant}[thm]{Variant}
\theoremstyle{definition}
\newtheorem{construction}[thm]{Construction}
\newtheorem{notations}[thm]{Notations}
\newtheorem{question}[thm]{Question}
\newtheorem{problem}[thm]{Problem}
\newtheorem{remark}[thm]{Remark}
\newtheorem{remarks}[thm]{Remarks}
\newtheorem{definition}[thm]{Definition}
\newtheorem{claim}[thm]{Claim}
\newtheorem{assumption}[thm]{Assumption}
\newtheorem{assumptions}[thm]{Assumptions}
\newtheorem{properties}[thm]{Properties}
\newtheorem{example}[thm]{Example}
\newtheorem{program}[thm]{Program}
\newtheorem{setup}[thm]{Set-up}
\newtheorem{conjecture}[thm]{Conjecture}
\catcode`\@=11
\def\opn#1#2{\def#1{\mathop{\kern0pt\fam0#2}\nolimits}}
\def\bold#1{{\bf #1}}%
\def\underrightarrow{\mathpalette\underrightarrow@}
\def\underrightarrow@#1#2{\vtop{\ialign{$##$\cr
 \hfil#1#2\hfil\cr\noalign{\nointerlineskip}%
 #1{-}\mkern-6mu\cleaders\hbox{$#1\mkern-2mu{-}\mkern-2mu$}\hfill
 \mkern-6mu{\to}\cr}}}
\let\underarrow\underrightarrow
\def\underleftarrow{\mathpalette\underleftarrow@}
\def\underleftarrow@#1#2{\vtop{\ialign{$##$\cr
 \hfil#1#2\hfil\cr\noalign{\nointerlineskip}#1{\leftarrow}\mkern-6mu
 \cleaders\hbox{$#1\mkern-2mu{-}\mkern-2mu$}\hfill
 \mkern-6mu{-}\cr}}}
\let\amp@rs@nd@\relax
\newdimen\ex@
\ex@.2326ex
\newdimen\bigaw@
\newdimen\minaw@
\minaw@16.08739\ex@
\newdimen\minCDaw@
\minCDaw@2.5pc
\newif\ifCD@
\def\minCDarrowwidth#1{\minCDaw@#1}
\newenvironment{CD}{\@CD}{\@endCD}
\def\@CD{\def\A##1A##2A{\llap{$\vcenter{\hbox
 {$\scriptstyle##1$}}$}\Big\uparrow\rlap{$\vcenter{\hbox{%
$\scriptstyle##2$}}$}&&}%
\def\V##1V##2V{\llap{$\vcenter{\hbox
 {$\scriptstyle##1$}}$}\Big\downarrow\rlap{$\vcenter{\hbox{%
$\scriptstyle##2$}}$}&&}%
\def\={&\hskip.5em\mathrel
 {\vbox{\hrule width\minCDaw@\vskip3\ex@\hrule width
 \minCDaw@}}\hskip.5em&}%
\def\verteq{\Big\Vert&&}%
\def\noarr{&&}%
\def\vspace##1{\noalign{\vskip##1\relax}}\relax\let\amp@rs@nd@&\iffalse}\fi
 \CD@true\vcenter\bgroup\relax\let\\=\cr\iffalse}\fi\tabskip\z@skip\baselineskip20\ex@
 \lineskip3\ex@\lineskiplimit3\ex@\halign\bgroup
 &\hfill$\m@th##$\hfill\cr}
\def\@endCD{\cr\egroup\egroup}
\def\>#1>#2>{\amp@rs@nd@\setbox\z@\hbox{$\scriptstyle
 \;{#1}\;\;$}\setbox\@ne\hbox{$\scriptstyle\;{#2}\;\;$}\setbox\tw@
 \hbox{$#2$}\ifCD@
 \global\bigaw@\minCDaw@\else\global\bigaw@\minaw@\fi
 \ifdim\wd\z@>\bigaw@\global\bigaw@\wd\z@\fi
 \ifdim\wd\@ne>\bigaw@\global\bigaw@\wd\@ne\fi
 \ifCD@\hskip.5em\fi
 \ifdim\wd\tw@>\z@
 \mathrel{\mathop{\hbox to\bigaw@{\rightarrowfill}}\limits^{#1}_{#2}}\else
 \mathrel{\mathop{\hbox to\bigaw@{\rightarrowfill}}\limits^{#1}}\fi
 \ifCD@\hskip.5em\fi\amp@rs@nd@}
\def\<#1<#2<{\amp@rs@nd@\setbox\z@\hbox{$\scriptstyle
 \;\;{#1}\;$}\setbox\@ne\hbox{$\scriptstyle\;\;{#2}\;$}\setbox\tw@
 \hbox{$#2$}\ifCD@
 \global\bigaw@\minCDaw@\else\global\bigaw@\minaw@\fi
 \ifdim\wd\z@>\bigaw@\global\bigaw@\wd\z@\fi
 \ifdim\wd\@ne>\bigaw@\global\bigaw@\wd\@ne\fi
 \ifCD@\hskip.5em\fi
 \ifdim\wd\tw@>\z@
 \mathrel{\mathop{\hbox to\bigaw@{\leftarrowfill}}\limits^{#1}_{#2}}\else
 \mathrel{\mathop{\hbox to\bigaw@{\leftarrowfill}}\limits^{#1}}\fi
 \ifCD@\hskip.5em\fi\amp@rs@nd@}
\newenvironment{CDS}{\@CDS}{\@endCDS}
\def\@CDS{\def\A##1A##2A{\llap{$\vcenter{\hbox
 {$\scriptstyle##1$}}$}\Big\uparrow\rlap{$\vcenter{\hbox{%
$\scriptstyle##2$}}$}&}%
\def\V##1V##2V{\llap{$\vcenter{\hbox
 {$\scriptstyle##1$}}$}\Big\downarrow\rlap{$\vcenter{\hbox{%
$\scriptstyle##2$}}$}&}%
\def\={&\hskip.5em\mathrel
 {\vbox{\hrule width\minCDaw@\vskip3\ex@\hrule width
 \minCDaw@}}\hskip.5em&}
\def\verteq{\Big\Vert&}
\def\novarr{&}
\def\noharr{&&}
\def\SE##1E##2E{\slantedarrow(0,18)(4,-3){##1}{##2}&}
\def\SW##1W##2W{\slantedarrow(24,18)(-4,-3){##1}{##2}&}
\def\NE##1E##2E{\slantedarrow(0,0)(4,3){##1}{##2}&}
\def\NW##1W##2W{\slantedarrow(24,0)(-4,3){##1}{##2}&}
\def\slantedarrow(##1)(##2)##3##4{%
\thinlines\unitlength1pt\lower 6.5pt\hbox{\begin{picture}(24,18)%
\put(##1){\vector(##2){24}}%
\put(0,8){$\scriptstyle##3$}%
\put(20,8){$\scriptstyle##4$}%
\end{picture}}}
\def\vspace##1{\noalign{\vskip##1\relax}}\relax\let\amp@rs@nd@&\iffalse}\fi
 \CD@true\vcenter\bgroup\relax\let\\=\cr\iffalse}\fi\tabskip\z@skip\baselineskip20\ex@
 \lineskip3\ex@\lineskiplimit3\ex@\halign\bgroup
 &\hfill$\m@th##$\hfill\cr}
\def\@endCDS{\cr\egroup\egroup}
\newdimen\TriCDarrw@
\newif\ifTriV@
\newenvironment{TriCDV}{\@TriCDV}{\@endTriCD}
\newenvironment{TriCDA}{\@TriCDA}{\@endTriCD}
\def\@TriCDV{\TriV@true\def\TriCDpos@{6}\@TriCD}
\def\@TriCDA{\TriV@false\def\TriCDpos@{10}\@TriCD}
\def\@TriCD#1#2#3#4#5#6{%
\setbox0\hbox{$\ifTriV@#6\else#1\fi$}
\TriCDarrw@=\wd0 \advance\TriCDarrw@ 24pt
\advance\TriCDarrw@ -1em
\def\SE##1E##2E{\slantedarrow(0,18)(2,-3){##1}{##2}&}
\def\SW##1W##2W{\slantedarrow(12,18)(-2,-3){##1}{##2}&}
\def\NE##1E##2E{\slantedarrow(0,0)(2,3){##1}{##2}&}
\def\NW##1W##2W{\slantedarrow(12,0)(-2,3){##1}{##2}&}
\def\slantedarrow(##1)(##2)##3##4{\thinlines\unitlength1pt
\lower 6.5pt\hbox{\begin{picture}(12,18)%
\put(##1){\vector(##2){12}}%
\put(-4,\TriCDpos@){$\scriptstyle##3$}%
\put(12,\TriCDpos@){$\scriptstyle##4$}%
\end{picture}}}
\def\={\mathrel {\vbox{\hrule
   width\TriCDarrw@\vskip3\ex@\hrule width
   \TriCDarrw@}}}
\def\>##1>>{\setbox\z@\hbox{$\scriptstyle
 \;{##1}\;\;$}\global\bigaw@\TriCDarrw@
 \ifdim\wd\z@>\bigaw@\global\bigaw@\wd\z@\fi
 \hskip.5em
 \mathrel{\mathop{\hbox to \TriCDarrw@
{\rightarrowfill}}\limits^{##1}}
 \hskip.5em}
\def\<##1<<{\setbox\z@\hbox{$\scriptstyle
 \;{##1}\;\;$}\global\bigaw@\TriCDarrw@
 \ifdim\wd\z@>\bigaw@\global\bigaw@\wd\z@\fi
 \mathrel{\mathop{\hbox to\bigaw@{\leftarrowfill}}\limits^{##1}}
 }
 \CD@true\vcenter\bgroup\relax\let\\=\cr\iffalse}\fi
 \tabskip\z@skip\baselineskip20\ex@
 \lineskip3\ex@\lineskiplimit3\ex@
 \ifTriV@
 \halign\bgroup
 &\hfill$\m@th##$\hfill\cr
#1&\multispan3\hfill$#2$\hfill&#3\\
&#4&#5\\
&&#6\cr\egroup%
\else
 \halign\bgroup
 &\hfill$\m@th##$\hfill\cr
&&#1\\%
&#2&#3\\
#4&\multispan3\hfill$#5$\hfill&#6\cr\egroup
\fi}
\def\@endTriCD{\egroup}
\newcommand{\sA}{{\mathcal A}}
\newcommand{\sB}{{\mathcal B}}
\newcommand{\sC}{{\mathcal C}}
\newcommand{\sD}{{\mathcal D}}
\newcommand{\sE}{{\mathcal E}}
\newcommand{\sF}{{\mathcal F}}
\newcommand{\sG}{{\mathcal G}}
\newcommand{\sH}{{\mathcal H}}
\newcommand{\sI}{{\mathcal I}}
\newcommand{\sJ}{{\mathcal J}}
\newcommand{\sK}{{\mathcal K}}
\newcommand{\sL}{{\mathcal L}}
\newcommand{\sM}{{\mathcal M}}
\newcommand{\sN}{{\mathcal N}}
\newcommand{\sO}{{\mathcal O}}
\newcommand{\sP}{{\mathcal P}}
\newcommand{\sQ}{{\mathcal Q}}
\newcommand{\sR}{{\mathcal R}}
\newcommand{\sS}{{\mathcal S}}
\newcommand{\sT}{{\mathcal T}}
\newcommand{\sU}{{\mathcal U}}
\newcommand{\sV}{{\mathcal V}}
\newcommand{\sW}{{\mathcal W}}
\newcommand{\sX}{{\mathcal X}}
\newcommand{\sY}{{\mathcal Y}}
\newcommand{\sZ}{{\mathcal Z}}
\newcommand{\A}{{\mathbb A}}
\newcommand{\B}{{\mathbb B}}
\newcommand{\C}{{\mathbb C}}
\newcommand{\D}{{\mathbb D}}
\newcommand{\E}{{\mathbb E}}
\newcommand{\F}{{\mathbb F}}
\newcommand{\G}{{\mathbb G}}
\newcommand{\HH}{{\mathbb H}}
\newcommand{\I}{{\mathbb I}}
\newcommand{\J}{{\mathbb J}}
\renewcommand{\L}{{\mathbb L}}
\newcommand{\M}{{\mathbb M}}
\newcommand{\N}{{\mathbb N}}
\newcommand{\BP}{{\mathbb P}}
\newcommand{\Q}{{\mathbb Q}}
\newcommand{\R}{{\mathbb R}}
\newcommand{\T}{{\mathbb T}}
\newcommand{\U}{{\mathbb U}}
\newcommand{\V}{{\mathbb V}}
\newcommand{\W}{{\mathbb W}}
\newcommand{\X}{{\mathbb X}}
\newcommand{\Y}{{\mathbb Y}}
\newcommand{\Z}{{\mathbb Z}}
\newcommand{\id}{{\rm id}}
\newcommand{\rk}{{\rm rank}}
\newcommand{\END}{{\mathbb E}{\rm nd}}
\newcommand{\End}{{\rm End}}
\newcommand{\Hg}{{\rm Hg}}
\newcommand{\tr}{{\rm tr}}
\newcommand{\Sl}{{\rm Sl}}
\newcommand{\Gl}{{\rm Gl}}
\newcommand{\Sp}{{\rm Sp}}
\newcommand{\MT}{{\rm MT}}
\newcommand{\Cor}{{\rm Cor}}
\newcommand{\Hom}{{\sH}{\rm om}}
\newcommand{\Mon}{{\rm Mon}}
\newcommand{\s}{{\rm sl}}
\newcommand{\ch}{{\rm c}}
\newcommand{\gr}{{\mathfrak g \mathfrak r}}
\maketitle
\section*{Introduction}\label{in}
Let $f:X\to Y$ be a family of $m$-folds over a curve $Y$, smooth over $U=Y\setminus S$.
We will assume that $f$ is semistable, i.e. that $X$ is non-singular and that $\Delta=f^{-1}(S)$ is a reduced
normal crossing divisor.

A slight improvement of the Arakelov inequality in \cite{VZ} has been shown in
\cite{MVZ}: 
\begin{theorem}\label{in.1}
For all $\nu\geq 1$ with $f_*\omega_{X/Y}^\nu\neq 0$ one has the inequality
\begin{equation}\label{eq.1}
\frac{\deg (f_* \omega^{\nu}_{X/Y})}{{\rm rank} (f_*
\omega^{\nu}_{X/Y})}\leq \frac{m\cdot \nu}{2} \cdot \deg(\Omega^1_Y(\log S)).
\end{equation}
\end{theorem}
Assume that the family $f$ is not isotrivial, hence that for no finite covering $Y'\to Y$ one has a birational map $X\times_YY' \to F\times Y'$. If one assumes in addition that the smooth fibres $F$ of $f$ are minimal models, 
then for $\nu >1$ the left hand side of (\ref{eq.1}) is strictly larger than zero. If $Y=\BP^1$ one sees
that $\#S\geq 3$.

Here we are interested in the case $\nu=1$, assuming that the non-isotriviality implies $\deg (f_* \omega_{X/Y}) >0$,
hence if a local Torelli Theorem holds for $F$. In this case we consider $f$
to be ``minimal'' if 
\begin{equation}\label{eq.2}
\frac{\deg (f_* \omega_{X/Y})}{{\rm rank} (f_*\omega_{X/Y})} = \frac{m}{2} \cdot \deg(\Omega^1_Y(\log S)).
\end{equation}
So a general question would be, whether there are such minimal families, and what is special about
them and their variation of Hodge structures. The case of families of curves (discussed in \cite{MVZ})
and the lack of any examples made us believe, that such minimal families only exists if
$\dim(H^0(F,\omega_F))=1$. Some examples of families of minimal families of Calabi-Yau manifolds have been
discussed in \cite{VZ1}.

Here we will restrict ourselves to families of $g$-dimensional Abelian varieties, and we will state some generalization
of (\ref{eq.1}) for families over a higher higher dimensional bases. Before defining the necessary tools,
let us give two examples, for $Y=\BP^1$. As well known, and as we will recall in Section \ref{ai}, for semistable families of Abelian varieties the minimal number of singular fibres is $\# S = 4$.    
\begin{example}\label{in.2} (Beauville \cite{Bea}) If $g:E\to \BP^1$ is a semistable non-isotrivial family
of elliptic curves with $4$ singular fibres, then it is a modular family. Moreover, up to isomorphism there
are $6$ possible examples.
\end{example}

\begin{example}\label{in.3}  (Theorem 0.2 in \cite {VZ2})
Let $f:X\to \mathbb P^1$ be a non-isotrivial semistable family of Abelian varieties of dimension $g$
with $4$ singular fibres. Then (\ref{eq.2}) says that $\deg (f_* \omega_{X/Y})=g$, and this 
implies that $f$ is isogenous to $E\times_{\BP^1} E 
\times_{\BP^1} \cdots \times_{\BP^1} E$ where $E\to {\BP^1}$ is one Beauville's examples.
\end{example}

One possible interpretation is that $\mathbb P^1\setminus \{y_1,\cdots,y_4\}$ is a ``minimal'' 
subvariety of the moduli stack $\sA_g$ of polarized Abelian variety, which forces the corresponding
families to be ``special''. 

Special families of Abelian varieties should become isogenous to a universal family over a Shimura subvariety of $\sA_g$, after replacing the base $U$ by some \'etale covering.
Here and afterwards ``Shimura variety'' stands for ``Shimura variety of
Hodge type'', as defined in \cite{Mum1} and \cite{Mum2}, i.e. for a moduli scheme
of Abelian varieties with a prescribed Mumford-Tate group.

In this note we want to speculate about the meaning of ``minimal'' for families over surfaces
or over higher dimensional bases, so we propose the following
\begin{program}\label{in.4} \ 
\begin{itemize}
\item Define minimal subvarieties $U$ in the moduli stack $\sA_g$ of polarized $g$-dimensional
Abelian varieties, in terms of numerical conditions of the corresponding
$\Q$-variation of Hodge structures.
\item Show that minimal subvarieties are Shimura varieties.
\item Try to distinguish different types of Shimura varieties by numerical
conditions. 
\end{itemize}
\end{program}
We are grateful to the referee for pointing out several ambiguities in the first version
of this article, and for suggestions how to improve the presentation.
\section{Variations of Hodge structures and Higgs bundles}\label{hh}
\begin{setup}\label{hh.1}
Let $Y$ denote an $n$-dimensional complex projective manifold, $S$ a
normal crossing divisor  and $U=Y\setminus S$. Let $\varphi: U\to \sA_g$ be a morphism
to the moduli stack of polarized Abelian $g$-dimensional varieties, and let $f:V\to U$ be the induced family of Abelian varieties (usually one can assume that $\varphi$ factors through a morphism to some fine moduli scheme 
and that $f$ is the pullback of the universal family).
Choose a non-singular compactification $X$ of $V$ such that $f$ extends to $f:X\to Y$,
and such that $T=f^{-1}(S)$ is a normal crossing divisor.

Let $R^1f_*\Q_V$ be the induced variation of Hodge structures on $U$. We will assume that
the local monodromy in general points of the components of $S$ is unipotent.
For example this will be true if $f:X\to Y$ is semistable in codimension one, i.e.
if $S$ contains a dense open subset $S^0$ such that $f^{-1}(S^0)$ is a reduced normal crossing divisor.
\end{setup}

Let $\V$ be a $\C$-subvariation of Hodge structures of $R^1f_*\C_V\otimes \C=R^1f_*\C_V$.
Consider the Deligne extension $\sH$ of the locally free sheaf
$\V\otimes_\Q\sO_U$ to $Y$, defined in \cite{Del}. This distinguished extension can be characterized  as follows:\\
Let  $U$ be an open subset in $Y$  with $U\cap S\not= \emptyset,$ $s$  a section
in $\sH_{U\setminus S}$  and $s_1,\cdots, s_n$ a multiple-valued flat frame, such that  
$s=\sum_if_is_i,$ with multiple-valued coefficient functions $f_i.$  Then $s$ extends holomorphically to $U$ if and only the $f_i$ have at most logarithmic singularities (\cite{Sch}, Page 235).
The sheaf $\sH$ is locally free and the Gauss-Manin connection extends to
$$\nabla:\sH \to \sH\otimes \Omega^1_Y(\log S).$$

As shown in \cite{Sch} the $F^\bullet$ filtration extends to a filtration of $\sH$ by subbundles.
Griffiths transversality implies that taking the graded sheaf 
$$F=\mathfrak{gr}_{F^\bullet}(\sH)=F^{1,0}\oplus F^{0,1}$$
one obtains a logarithmic Higgs bundle, i.e. a locally free sheaf $F$ together
with the Higgs field $\tau: F\to F\otimes \Omega^1_Y(\log S)$,
an $\sO_Y$ linear map induced by the Gauss-Manin connection. 
Here $\tau|_{F^{0,1}}$ is zero, and we will write
$$\tau: F^{1,0}\>>> F^{0,1}\otimes \Omega^1_Y(\log S).$$
Note that $\tau$ is just the cup product with the Kodaira-Spencer class.

For a family of Abelian varieties, $R^rf_*\C_V$ is the $g$-th wedge product of
$R^1F_*\C_V$. Using the notations introduced above for $\V=R^1f_*\C_V$, the iterated
cup product with the Kodaira-Spencer map gives 
\begin{multline*}
\tau^\ell:\bigwedge^g F^{1,0}\>>> \bigwedge^{g-1}F^{1,0}\otimes F^{0,1} \otimes\Omega^1_Y(\log S)
\>>> 
\\
\cdots \>>> \bigwedge^{g-\ell} F^{1,0}\otimes \bigwedge^{\ell} F^{0,1}\otimes S^{\ell}\Omega^1_Y(\log S).
\end{multline*}
For $\ell=g$ one obtains the Griffiths Yukawa coupling
$$ 
\tau^g:\bigwedge^g F^{1,0}\>>> \bigwedge^g F^{0,1}\otimes S^{g}\Omega^1_Y(\log S).
$$
Candidates for the numerical invariants we are looking for in Program \ref{in.4}
could be $\deg(F^{1,0})$, $\deg(\Omega^1_Y(\log S))$ and the maximal number $\ell$ with
$\tau^\ell\neq 0$. 

\section{Arakelov inequalities over curves}\label{ai}
Let us return to families over a curve $Y$ in Set-up \ref{hh.1}, and to a $\C$-subvariation
of Hodge structures $\V\subset R^1f_*\C_V$ with Higgs bundle $(F,\tau)$. The inequality
(\ref{eq.1}) for $\nu=1$ is a special case of Faltings Arakelov inequality:
\begin{theorem}[\cite{Fal} for $\V=R^1f_*\C_V$ and \cite{Del}, in general] \label{ai.1} \ \\
If $\dim U=1$ one has the inequality
\begin{equation}\label{eq.3}
2\cdot \deg(F^{1,0})\leq \rk(F^{1,0}) \cdot \deg(\Omega^1_Y(\log S)).
\end{equation}
\end{theorem}
Assume for a moment that $\T$ is the largest unitary local subsystem in $R^1f_*\C_V$
and $\V\oplus \T=R^1f_*\C_V$. Simpson's correspondence \cite{S90} implies that
$F^{1,0}$ can not have an invertible quotient of degree $\leq 0$. Hence
for $Y=\BP^1$ the sheaf $F^{1,0}={F^{0,1}}^\vee$ has to be the direct sum
of line bundles $\sO_Y(\nu)$ for $\nu>0$. Obviously this implies that
$\deg(\Omega^1_Y(\log S))\geq 2$. Simpson's correspondence is also the main ingredient
in the proof of the next result:
\begin{theorem}[\cite{VZ2}] \label{ai.2} In Theorem \ref{ai.1}
assume that
\begin{equation}\label{eq.4}
2\cdot \deg(F^{1,0})= g \cdot \deg(\Omega^1_Y(\log S)).
\end{equation}
Then $F^{1,0}$ and $F^{0,1}$ are polystable and
$\tau: F^{1,0}\to F^{0,1}\otimes \Omega^1_Y(\log S)$ is an isomorphism.
\end{theorem}
\begin{theorem}[\cite{VZ2}]\label{ai.3}
Assume that $\V=R^1f_*\C_V$. If for the Higgs field $(F,\tau)$ of $\V=R^1f_*\C_V$ (\ref{eq.4})
holds, then $U$ is a rigid Shimura curve. Replacing $U$ by an \'etale covering, $V\to U$ is isogenous to
a selfproduct of the corresponding universal family. 
\end{theorem}
Here rigid means that there is no non-trivial deformation of the morphism from $U$ to the moduli stack
$\mathcal A_g$. This type of Shimura curves has first been constructed by Mumford \cite{Mum2}  for $g=4$.

Of course, if $\tau: F^{1,0}\to F^{0,1}\otimes \Omega^1_Y(\log S)$ is an isomorphism,
(\ref{eq.3}) is an equality. In \cite{VZ2} the polystability was expressed in a different way.
If $F^{1,0}$ is polystable, replacing $U$ by an \'etale covering, one may write
$$
F^{1,0}=\sL\otimes \sU \mbox{ \ \ and \ \ } F^{0,1}=\sL^{-1}\otimes \sU^\vee
$$
where $\sL$ is a logarithmic theta characteristic, i.e. an invertible sheaf
with $\sL^2=\omega_Y(S)$, and where $\sU$ is a unitary locally free sheaf.

The Example \ref{in.3} follows immediately from:
\begin{addendum}[\cite{VZ2}]\label{ai.4}
If $S\neq \emptyset$ then the conclusion of Theorem \ref{ai.3} can be strengthened:\\[.1cm]
$U$ is a modular curve and (replacing $U$ by an \'etale covering) $V\to U$ is isogenous to the $g$-fold selfproduct of the corresponding modular family of elliptic curves.
\end{addendum}
For $S=\emptyset$ the description of the isogeny class of $V\to U$ is more complicated.
Proving \ref{ai.2} in \cite{VZ2}, one recovers the description of the variation of
Hodge structures $R^1f_*\C_V$ as the corestriction of a quaternion algebra defined over a totally real
number field which is ramified at all infinite places except one.

Remark that in \cite{VZ2} we also consider the non-rigid case. Here one has to decompose
$R^1f_*\C_V$ as $\T\oplus\V$, with $\T$ unitary. If for the logarithmic Higgs field of
$\V$ the equality (\ref{eq.4}) holds $U$ is a Shimura curve and the family $V\to U$ is isogenous to the product of certain universal families. We skip the details, since for $\dim(U)>1$
we are at the moment not able to allow unitary direct factors in the variation of Hodge structures.

\section{Arakelov inequalities over a higher dimensional base}\label{ar}
Let $Y$ be an $n$-dimensional projective manifold.
Define for a torsion free coherent sheaf $\sF$ on $Y$
$$
\Upsilon(\sF)=\frac{\ch_1(\sF)}{\rk(\sF)} \in H^2(Y,\Q),
$$
and the discriminant
$$
\Delta(\sF)= 2\cdot \rk(\sF) \cdot \ch_2(\sF) - (\rk(\sF)-1)\cdot \ch_1(\sF)^2 \in H^4(Y,\Q).
$$
Choose on $Y$ an invertible sheaf $\sN$, or more generally
an $\R$-divisor $\sN$. One defines the slope and the discriminant of a torsion free coherent sheaf $\sF$ as
$$
\mu_\sN(\sF)=\Upsilon(\sF).\ch_1(\sN)^{n-1} \mbox{ \ \ and \ \ }
\delta_\sN(\sF)=\Delta(\sF).\ch(\sN)^{n-2},
$$
respectively. Both, $\Upsilon$ and $\mu_\sN$, are additive for tensor products.
\begin{definition}\label{ar.1} Let $\sH$ be an ample invertible sheaf on $Y$.
\begin{enumerate}
\item[i.] $\sF$ is numerically effective (nef), if for all curves $\tau:C\to Y$ and for all
invertible quotient sheaves $\sL$ of $\tau^*\sF$ one has $\deg(\sL)\geq 0$.
\item[ii.] $\sF$ is ample with respect to an open subscheme $U'$ of $Y$, if for some $\nu \gg 0$
there exists a morphism $\oplus \sH \to S^\nu(\sF)$, which is an surjective over $U'$.

If $\sF$ is invertible, this is equivalent to:\\
For some $\mu > 0$ the sheaf $\sF^\mu$ is generated by $H^0(Y,\sF^\mu)$ in all points $u\in U'$,
and the induced morphism $U'\to \BP(H^0(Y,\sF^\mu)$ is an embedding.
\item[iii.] $\sF$ is big, if it is ample with respect to some open dense subscheme.
\end{enumerate}
\end{definition}
Let us also recall the notion of stability for a torsion free coherent sheaf $\sF$ with respect to the
invertible sheaf or an $\R$-divisor $\sN$.
\begin{enumerate}
\item[iv.] $\sF$ is $\mu_\sN$-stable, if for all subsheaves $\sG$ with $\rk(\sG)<\rk(\sF)$ one has
$$
\mu_\sN(\sG) < \mu_\sN(\sF).
$$
\item[v.] $\sF$ is $\mu_\sN$-semistable, if for all subsheaves $\sG$ one has
$$
\mu_\sN(\sG) \leq \mu_\sN(\sF).
$$
\item[vi.] A $\mu_\sN$-semistable sheaf $\sF$ is $\mu_\sN$-polystable, if it is the direct sum of
$\mu_\sN$-stable sheaves.
\end{enumerate}
In iv), v) and vi) one should assume at least that $\sN$ is nef and big.\par

Consider the following set of assumptions for $(Y,S)$:
\begin{assumption}\label{ar.2} \
\begin{enumerate}
\item[a.] $\Omega^1_Y(\log S)$ is nef.
\item[b.] $\omega_Y(S)$ is ample with respect to $U$.
\item[c.] $\omega_Y(S)$ is semiample.
\item[d.] If $\sF$ and $\sG$ are two $\mu_{\omega_Y(S)}$ polystable locally free sheaves,
then $\sF\otimes \sG$ is again $\mu_{\omega_Y(S)}$ polystable.
\item[e.] If $\sF$ is $\mu_{\omega_Y(S)}$ polystable  and locally free,
then
$$
\mu_{\omega_Y(S)}(\sF)=\delta_{\omega_Y(S)}(\sF)=0
$$
if and only if $\sF$ is unitary.
\end{enumerate}
\end{assumption}
\begin{remark}\label{ar.3}
The condition \ref{ar.2} a) implies that $\omega_Y(S)$ is nef, so a) and b) together imply that
$\omega_Y(S)$ nef and ample with respect to $U$. The conditions d) and e) hold true whenever $\omega_Y(\log S)$ is
ample. As Yau told us, he and Sun are able to show that a) and b) imply d) and e).

The additional condition c) implies that Simpson's correspondence remains true if one chooses
$\omega_Y(S)$ as polarization:
\begin{enumerate}
\item[f.] Let $(E\theta)$ be a logarithmic Higgs bundle induced by a $\C$ variation of Hodge structures $\V$
on $U$ with unipotent local monodromy. If $G\subset E$ is a saturated Higgs subsheaf then
\begin{equation}\label{eq.5}
\ch_1(G).\ch_1(\omega_Y(S))^{n-1}\leq 0,
\end{equation}
and the equality holds, if and only if $G$ is induced by a local subsystem.
\end{enumerate}
Here ``saturated'' means that the quotient sheaf has no torsion, and $\sG$ is a Higgs
subsheaf, if $\theta(\sG)\subset \sG\otimes \Omega^1_Y(\log S)$. If $\sG$ is induced by a local subsystem of $\V$ then $\sG$ is of course a 
Higgs subbundle. 

In fact consider for $r=\rk G$ the rank one Higgs subbundle
$$
(\det G,0)=\bigwedge^r(G,\theta) \subset \bigwedge^r(E,\theta).
$$
The curvature of the Hodge metric $h$ on $\det G$ is negative semidefinite, and
the Chern form $\ch_1(G,h)$ represents the Chern class of $\det (G)$.
One obtains (\ref{eq.5}).

Assume now, that (\ref{eq.5}) is an equality.
If $\ch_1(G,h)\neq 0$ one finds an irreducible non-singular curve $C$ as the zero locus of $n-1$ general
sections of some power of $\omega_Y(S)$, such that $C$ does not lie in the kernel of
$\ch_1(G,h)$. So the pullback of $\ch_1(G,h)$ to $C$ is non-zero. This implies that
at least in one point of $C$ the pullback of $\ch_1(G,h)$ is strictly negative, hence
$\ch_1(G)\cdot \ch_1(\omega_Y(S))^{n-1}<0$, contrary to the assumption.
So $\ch_1(G,h)=0$ hence $\ch_1(G).\ch_1(\sH)^{n-1}=0$ for all ample invertible sheaves $\sH$,
and \cite{Sim} implies that $G$ is induced by a local system.
\end{remark}

Since $\omega_Y(S)$ is nef and ample with respect to $U$ (assuming \ref{ar.2}, a) and b)),
one can apply Yau's uniformization theorem (\cite{Yau2} and \cite{Yau}), recalled in \cite[Section 1]{VZ4}. In particular, the sheaf $\Omega_Y^1(\log S)$ is $\mu_{\omega_Y(S)}$-polystable. Hence one has a direct sum decomposition
\begin{equation}\label{eq.6}
\Omega_Y^1(\log S)=\Omega_1\oplus \cdots \oplus \Omega_{s''}\oplus \cdots \oplus \Omega_{s'}
\oplus \cdots \oplus \Omega_s
\end{equation}
in stable sheaves.

\begin{notations}\label{ar.4} We choose the indices such that
\begin{enumerate}
\item[I.] for $i=1, \ldots , s''$ the sheaf $\Omega_i$ is invertible.
\item[II.] for $i=s''+1, \ldots , s'$, the sheaf $\Omega_i$ is stable
of rank $>1$ as well as $S^m(\Omega_i)$, for all $m>1$.
\item[III.] for $i=s'+1, \ldots , s$ we have the remaining stable direct factors,
i.e. those with $S^{m_i}(\Omega_i)$ non-stable, for some $m_i>1$.
\end{enumerate}
\end{notations}

Let $\V$ be a $\C$-subvariation of Hodge structures of $\W$, without a unitary direct factor, and with Higgs bundle
$$
\big(E^{1,0}\oplus E^{0,1} \ , \ \theta: E^{1,0}\to E^{0,1}\otimes \Omega^1_Y(\log S)\big).
$$
By \cite[Theorem 1]{VZ4}, assuming \ref{ar.2}, one has the Arakelov type inequality
\begin{equation}\label{eq.7}
\mu_{\omega_Y(S)}(E^{1,0})-\mu_{\omega_Y(S)}(E^{0,1}) \leq \mu_{\omega_Y(S)}(\Omega^1_Y(\log S)).
\end{equation}
Obviously, for $\dim(Y)=1$ this is the same as the inequality (\ref{eq.3}).\par

Unfortunately at present only part of Theorem \ref{ai.2} generalizes to higher dimensional base schemes
satisfying \ref{ar.2}, and to irreducible subvariations of Hodge structures $\V$.

By \cite[Theorem 1]{VZ4} the Arakelov equality
\begin{equation}\label{eq.8}
\mu_{\omega_Y(S)}(E^{1,0})-\mu_{\omega_Y(S)}(E^{0,1}) = \mu_{\omega_Y(S)}(\Omega^1_Y(\log S))
\end{equation}
implies that $E^{1,0}$ and $E^{0,1}$ are both $\mu_{\omega_Y(S)}$-semistable.
In different terms, if the Higgs field
$$
\theta: E^{1,0} \>>> E^{0,1}\otimes \Omega^1_Y(\log S)
$$
is a morphism between locally free sheaves of the same $\mu_{\omega_Y(S)}$-slope, then
$E^{1,0}$, $E^{0,1}$, and hence $E^{0,1}\otimes \Omega^1_Y(\log S)$, are all $\mu_{\omega_Y(S)}$-semistable.
We expect more.

\begin{conjecture}\label{ar.5} Assuming \ref{ar.2}, let $\V$ be an irreducible subvariation of Hodge structures satisfying the Arakelov equality (\ref{eq.8}). 
\begin{enumerate} 
\item[a.] There exists a unique $i\in \{1, \ldots , s\}$ such that the Higgs field $\theta$ factors through
$$
\theta: E^{1,0}\>>> E^{0,1}\otimes \Omega_i\>>> E^{0,1}\otimes \Omega^1_Y(\log S).
$$
\item[b.] $E^{1,0}$ and $E^{0,1}$ are $\mu_{\omega_Y(S)}$-stable.
\end{enumerate}
\end{conjecture}
Recall that for a $\mu_{\omega_Y(S)}$-semistable locally free sheaf $\sF$ one has the Bogomolov inequality
$\delta_{\omega_Y(S)}(\sF)\geq 0.$ Then \ref{ar.2} and the Arakelov equality (\ref{eq.8}) imply that
$$
\delta_{\omega_Y(S)}(E^{1,0}) \geq 0 \mbox{ \ \ and \ \ } \delta_{\omega_Y(S)}(E^{0,1}) \geq 0.
$$
We say that $\V$ or $(E^{1,0}\oplus E^{0,1}, \theta)$ satisfy the Bogomolov equality,
if
\begin{equation}\label{eq.9}
{\rm Min}\{\delta_{\omega_Y(S)}(E^{1,0}),\delta_{\omega_Y(S)}(E^{0,1})\}= 0.
\end{equation}
The conjecture \ref{ar.5} has been verified in \cite[Proposition 4]{VZ4} under the additional assumption that all
the direct factors of $\Omega^1_Y(\log S)$ are of type I, or that there are no direct factors of type III,
and that the Higgs bundles of all direct factors $\V$ of $R^1f_*\C_V$ satisfy the Bogomolov equality. In fact, this is an intermediate step in the proof of:
\begin{theorem}[ {\cite[Theorem 6]{VZ4}} ]\label{ar.6}
In the Set-up \ref{hh.1} assume that \ref{ar.2} holds, that the induced morphism $U\to \sA_g$ to the moduli stack of
polarized Abelian varieties is generically finite, and that $s=s'$, i.e. that there are no direct factors
of $\Omega_i$ of $\Omega^1_Y(\log S)$ of type III.
Assume moreover that for all irreducible $\C$-subvariations of Hodge structures $\V$ of $R^1f_*\C_V$ with logarithmic Higgs bundle $(E^{1,0} \oplus E^{0,1},\theta)$ the Arakelov equality (\ref{eq.8}) holds
as well as one of the following conditions
\begin{enumerate}
\item[i.] ${\rm Min}\{ \delta_{\omega_Y(S)}(E^{1,0}), \delta_{\omega_Y(S)}(E^{0,1})\}=0,$
\item[ii.] $s''=s$, hence $\Omega_i$ is invertible for $i=1,\ldots,s$.
\end{enumerate}
Then $U$ is a rigid Shimura subvariety of the moduli stack $\sA_g$, and
the universal covering $\tilde U$ is the product of $s$ complex balls of dimensions $n_j=\rk(\Omega_j)$.
\end{theorem}
The proof of Theorem \ref{ar.6} is based on Yau's Uniformization Theorem, the Yau inequality
for direct factors of $\Omega^1_Y(\log S)$ of type II, on the Simpson correspondence
(see \cite{Sim} and \cite{Sim2}), and on certain characterizations of Shimura varieties,
due to Mumford (see \cite{Mum1}, \cite{Mum2}, and also \cite{Moo}). As a byproduct one
obtains a classification of the corresponding Higgs bundles, hence of the local systems $\V$.
Roughly speaking, they are the standard ones, tensorized with unitary local systems (see Section \ref{su} for the 2-dimensional case).

An affirmative answer to the conjecture \ref{ar.5} would allow to describe
the irreducible subvariations of Hodge structures $\V$ of $R^1f_*\C_V$
according to the type of the direct factor $\Omega_i$. 
Note that the decomposition in (\ref{eq.6}) corresponds to a decomposition
$$
\tilde{U}=M_1\times \cdots \times M_s
$$
of the universal covering $\tilde{U}$ of $U$.
If $\Omega_i$ is invertible, $M_i$ is a one dimensional complex ball.
If $\Omega_i$ has rank $n_i >1$ and if $S^m(\Omega_i)$ is stable for all $m>0$,
one can show that the existence of the variation of Hodge structures $\V_i$
in \ref{ar.5} with Arakelov equality (\ref{eq.8}) and with the Bogomolov
equality (\ref{eq.9}) implies that $M_i$ is a complex ball of dimension $n_i$. In both cases the results in \cite{VZ4} give an explicite description of $\V$. 
If $\Omega_i$ has rank $n_i >1$ and if $S^m(\Omega_i)$ is not stable for some $m>1$ the factor $M_i$ is a bounded symmetric domain of rank $>1$ and the Margulis Superrigidity holds. So one has additional tools to understand this case. 

\section{Characterization of generalized Hilbert modular surfaces, and Picard modular surfaces}\label{su}
Let us restrict ourselves to the case $\dim(Y)=2$. There are no direct factors
of type III. In fact, by \cite{Yau} if $S^m(\Omega^1_Y(\log S))$ is not stable, for some $m>0$, then $\Omega^1_Y(\log S)= \Omega_1\oplus \Omega_2$ with $\Omega_i$ invertible. So there are only two possible types of Shimura varieties:
\begin{example}\label{su.1} \ \\[.1cm]
{Type I.:} \ $U$ is an arithmetic quotient of the product of two copies of the upper half
plane, and $Y$ is a smooth toroidal compactification. So
$$
\Omega^1_Y(\log S)=\Omega_1\oplus\Omega_2,\mbox{ \ \ with \ \ } \mu_{\omega_Y(S)}(\Omega_1)=\mu_{\omega_Y(S)}(\Omega_2).
$$
Assume that (replacing $U$ by an \'etale covering, if necessary) there exists a variation $\L_i$ of Hodge structures on $U$ with Higgs bundle
$$
\Omega_i^{1/2} \ \oplus \ \Omega_i^{-1/2},\mbox{ \ \ and \ \ }
\theta_i: \Omega_i^{1/2} \> {\rm id} >> \Omega_i^{-1/2}\otimes
\Omega_i \subset \Omega_i^{-1/2}\otimes \Omega^1_Y(\log S).
$$
We will call $U$ a generalized Hilbert modular surface, and $\L_i$ the standard uniformizing variation of Hodge structures.

One easily verifies that the equality (\ref{eq.8}) holds for $\L_i$, and also for $\L_1\oplus\L_2$.\\[.2cm]
{Type II.:} \ $U$ is an arithmetic quotient of a two dimensional complex ball. Then
$S^m(\Omega^1_Y(\log S))$ is $\mu_{\omega_Y(S)}$ stable for all $m$. Assume there exists
a variation $\L$ of Hodge structures on $U$ with Higgs field
$$
\omega_Y(S)^{-1/3}\otimes \Omega^1_Y(\log S) \ \oplus \ \omega_Y(S)^{-1/3},
$$
again after replacing $U$ by an \'etale covering. Here the Higgs field is the identity
$$
\theta:\omega_Y(S)^{-1/3}\otimes \Omega^1_Y(\log S) \> {\rm id}>> \omega_Y(S)^{-1/3} \otimes \Omega^1_Y(\log S).
$$
We will call $U$ a Picard modular surface, and and $\L$ the standard uniformizing variation of Hodge structures.

Obviously the Arakelov equality (\ref{eq.8}) holds for $\L$ and for $\L^\vee$, as well as
the Bogomolov equality (\ref{eq.9}).
However (\ref{eq.7}) is a strict inequality for the real variation of Hodge structures $\L\oplus\L^\vee$.
\end{example}
Note that part a) of conjecture \ref{ar.5} holds true for surfaces, i.e. the factorization of the Higgs field. In fact, the only non-trivial case is the one of a generalized Hilbert modular surface,
which has been handled in \cite[Proposition 3.4]{VZ4} without using the 
Bogomolov equality (\ref{eq.8}). For Hilbert modular surfaces, part b) has been shown as well.
\begin{theorem}\label{su.2}
Let $f:X\to Y$ be a semistable family of
$g$-dimensional Abelian varieties over an algebraic surface $Y$, smooth over the complement of a normal crossing divisor
$S$. Assume that $(Y,S)$ satisfies the assumption \ref{ar.2}, and that the induced morphism $U\to \sA_g$ to the moduli stack of polarized Abelian varieties is generically finite. Assume moreover that (\ref{eq.8}) and (\ref{eq.9})
hold for all $\C$ irreducible subvariations of Hodge structures $\V$ in $R^1f_*\C_V$.

Then $U$ is a rigid Shimura variety with universal family $V\to U$ and (replacing $U$ by an \'etale covering) either
\begin{enumerate}
\item[a.] $U$ is a generalized Hilbert modular surface and
$$
R^1f_*\C_V= \L_1\otimes \U_1 \oplus \L_2\otimes \U_2
$$
for unitary local systems $\U_i$ and for the standard uniformizing variations of Hodge structures
$\L_i$, or
\item[b.] $U$ is a ball quotient and
$$
R^1f_*\C_V= \L\otimes \U \oplus \L^\vee \otimes \U^\vee
$$
where $\U$ is a unitary local system and $\L$ the standard uniformizing variations of Hodge structures.
\end{enumerate}
\end{theorem}
Theorem \ref{su.2} is just a special case of Theorem \ref{ar.6}, together with the description of the
variation of Hodge structures, given in \cite[Theorem 5]{VZ4}.

At present we only know the assumptions \ref{ar.2}, d) and e), to hold when $\omega_Y(S)$
is ample. This assumption is only reasonable if $S=\emptyset$. On the other hand,
$\dim(Y)=2$ implies that $\delta_{\omega_Y(S)}(E^{1,0})=\Delta(E^{1,0})$ is independent of the
polarization. Let $H$ be an ample divisor on $Y$ for which $\Omega_Y^1(\log S)$ 
is $\mu_H$ polystable. As shown in \cite[Section 8]{VZ4}, such an ample sheaf always
exists. If one assumes the equation (\ref{eq.8}) to hold for the slope with respect to
$\ch_1(\omega_Y(S))+\epsilon\cdot H$, for some $\epsilon_0>0$ and all $\epsilon_0 \geq \epsilon \geq 0$, then the conclusion of Theorem \ref{su.2} remains true, as shown in \cite[Variant 8.4]{VZ4}.

\begin{remark}\label{su.3}
Since we allow ourselves in Theorem \ref{su.2} to replace $U$ by an \'etale covering, we may as well assume that
unitary local systems $\U_i$ or $\U$ are given by a unitary representation of the fundamental group of $Y$.
In fact, if $S\neq\emptyset$ and if no \'etale covering of $U$ is the product of two curves,
the representation of $\pi_1(U)$ defining $\U$ has finite image.

If $U=U_1\times U_2$ one is in case a), and $\U_i$ is given by a representation of
$\pi_1(U_i)$. Again, either $U_i$ is compact, or the bundle trivializes over an \'etale covering.
\end{remark}
\section{The Griffiths-Yukawa coupling}\label{gy}
As remarked in Example \ref{su.1} II) the Arakelov equality will never hold for
the whole variation of Hodge structures given by a
family of $g$-dimensional Abelian varieties over Picard modular surfaces, but just for the irreducible
$\C$-subvariations. Moreover it is easy to see that the Griffiths-Yukawa coupling $\tau^g$ on the middle cohomology
$\wedge^gR^1f_*\C_V=R^gf_*\C_V$ has to vanish in this case. This can be exploited to give another characterization of Hilbert modular surfaces, using the degree of the invertible sheaf
$$
\bigwedge^g(F^{g,0})=\bigwedge^g f_*\Omega_Y^1(\log S)= f_*\omega_{X/Y},
$$
and the Griffiths-Yukawa coupling introduced in Section \ref{hh}, but not the discriminant $\Delta(F^{1,0})$. Recall that the Arakelov inequality for the Higgs bundle
$(F,\tau)$ of $R^1f_*\C_V$ says
\begin{equation}\label{eq.10}
\ch_1(f_*\omega_{X/Y}). \ch_1(\omega_Y(S))\leq \frac{g}{4}\ch_1(\omega_Y(S))^2.
\end{equation}
In fact one may write $R^1f_*\C_V=\V\oplus\T$ where $\T$ is the maximal unitary sublocal system.
If $(E,\theta)$ denotes the Higgs bundle of $\V$ then 
$$
\ch_1(E^{1,0}).\ch_1(\omega_Y(S))=\ch_1(F^{1,0}).\ch_1(\omega_Y(S))=\ch_1(f_*\omega_{X/Y}). \ch_1(\omega_Y(S))
$$ 
and (\ref{eq.10}) is nothing but (\ref{eq.7}). If $\T\neq 0$ the unequality (\ref{eq.10})
has to be strict.
\begin{proposition}\label{gy.1}
Let $f:X\to Y$ be a semistable family of $g$-dimensional Abelian varieties over an algebraic surface $Y$, smooth over the complement of a normal crossing divisor $S$. Assume that $(Y,S)$ satisfies
the assumption \ref{ar.2}, and that the induced morphism $U\to \sA_g$ to the moduli stack of
polarized Abelian varieties is generically finite. Then the following conditions are equivalent:
\begin{enumerate}
\item[i.] $U$ is a generalized Hilbert modular surfaces. Replacing $U$ by an \'etale covering and using the notations from \ref{su.1}
$$
R^1f_*\C_V= \L_1\otimes \U_1 \oplus \L_2\otimes \U_2
$$
for unitary local systems $\U_i$.\vspace{.1cm}
\item[ii.] $\displaystyle\ch_1(f_*\omega_{X/Y}). \ch_1(\omega_Y(S)) = \frac{g}{4}\ch_1(\omega_Y(S))^2$ and $\delta(F^{1,0})=0$. \vspace{.1cm}
\item[iii.]  $\displaystyle\ch_1(f_*\omega_{X/Y}). \ch_1(\omega_Y(S)) = \frac{g}{4}\ch_1(\omega_Y(S))^2$ and $\tau^g\neq 0$.
\end{enumerate}
\end{proposition}
\begin{proof}
It is obvious that i) implies ii) and iii). In \cite[Corollary 8.1]{VZ4}
it is shown that iii) implies i). 

Recall that the Griffiths-Yukawa coupling $\tau^g$ is a morphism 
$$
\tau^g:\det(F^{1,0}) \>>> \det(F^{0,1})\otimes S^g(\Omega^1_Y(\log S))
$$
between polystable sheaves of the same slope. If $\tau^g$ is non-zero,
$S^g(\Omega^1_Y(\log S))$ will be non-stable, hence $\Omega^1_Y(\log S)$ must be 
the direct sum of two invertible sheaves. 
\end{proof}

A characterization of Picard modular surfaces, without using the discriminant
and without the assumption that the Arakelov equality (\ref{eq.8}) holds for all irreducible $\C$-subvariations of Hodge structures of $R^1f_*\C_V$ is more difficult (and presumably not possible for $g$ large). Before stating a criterion for the lowest possible value $g=3$, let us introduce some 
notation. 
$$
(G=\bigoplus_{q=0}^g G^{g-q,q},\tau'|_G)
$$
denotes the Higgs subbundle of
$$
(F=\bigwedge^g F^{1,0}\oplus \bigwedge^{g-1}F^{1,0}\otimes F^{0,1} \oplus \cdots \oplus F^{1,0}\otimes \bigwedge^{g-1} F^{0,1} \oplus \bigwedge^g F^{0,1}, \tau')
$$
generated by $\det(F^{1,0})=\bigwedge^g F^{1,0}$. 
\begin{proposition}\label{gy.2}
Assume in \ref{gy.1} that $g=3$ and that $\tau^3=0$.
Then 
$$\ch_1(f_*\omega_{X/Y}). \ch_1(\omega_Y(S)) \leq \frac{2}{3}\ch_1(\omega_Y(S))^2
$$
and the following two conditions are equivalent:
\begin{enumerate}
\item[i.] $U$ is a Picard modular surface. Replacing $U$ by an \'etale covering and using again the notations from \ref{su.1}
$$
R^1f_*\C_V= \L\otimes \U \oplus \L^\vee \otimes \U^\vee
$$
for a unitary local system $\U$.
\item[ii.] \ \ \hspace*{\fill} $\displaystyle\ch_1(f_*\omega_{X/Y}). \ch_1(\omega_Y(S)) = \frac{2}{3}\ch_1(\omega_Y(S))^2.$
\hspace*{\fill} \ \
\end{enumerate}
\end{proposition}
\begin{proof}
Since $\varphi: U \to \sA_g$ is generically finite,
$$
\det(F^{1,0})\otimes T_Y(-\log S) \>>> \bigwedge^{2}F^{1,0}\otimes F^{0,1}
$$
must be injective. This implies that $G^{2,1}=\det(F^{1,0})\otimes T_Y(-\log S)$.
Since $\tau^3=0$ one has $G^{0,3}=0$ and 
$$
G=\det(F^{1,0})\oplus \det(F^{1,0})\otimes T_Y(-\log S) \oplus G^{1,2}
$$
is a Higgs subbundle. One finds
$$
\mu(\det(F^{1,0})) + 2\cdot \mu(\det(F^{1,0})) - 2 \cdot \mu(\Omega^1_Y(\log S)) +
\rk(G^{1,2})\cdot \mu(G^{1,2})\leq 0.
$$
$(G^{1,2},0)$ is a Higgs subbundle of the Higgs bundle of the local system $R^3f_*\C_V$, but it can not split. Hence $\mu(G^{1,2}) < 0$.
Since $\det(F^{1,0})\otimes S^2(T_Y(-\log S))$ is polystable, and since $G^{1,2}$
is a quotient sheaf,
$$
\mu(\det(F^{1,0})) - 2 \cdot \mu(\Omega_Y^1(\log S)) \leq \mu(G^{2,1}).
$$
Alltogether one finds
\begin{multline*}
6\cdot \big(\ch_1(f_*\omega_{X/Y}). \ch_1(\omega_Y(S)) - \frac{2}{3}\ch_1(\omega_Y(S))^2\big)=
6\cdot \mu(\det(F^{1,0})) - 8\cdot\mu(\Omega^1_Y(\log S))=\\
\mu(\det(F^{1,0})) + 2\cdot \mu(\det(F^{1,0})) - 2 \cdot \mu(\Omega^1_Y(\log S)) +\\
3\cdot\mu(\det(F^{1,0})) - 6 \cdot \mu(\Omega_Y^1(\log S))  \leq 0,
\end{multline*}
hence the inequality stated in d). Assume now that ii) holds, hence that
this is an equality. Going backwards, one sees that this is only possible if $\rk(G^{1,2})=3$,
hence if
$$
\det(F^{1,0})\otimes S^2(T_Y(\log S))=G^{1,2}.
$$
Moreover, by \ref{ar.3} $G$, as a Higgs subbundle of degree zero, is induced by some local system. Therefore $\ch_1(G)=0$ and  $\Delta(G)=0$. The discriminant does not change by tensorization with invertible sheaves, hence
$$
\Delta(G) = \Delta(\det(F^{1,0})\otimes S^2(\sO_Y\oplus \Omega^1_Y(\log S)))=
\Delta(S^2(\sO_Y\oplus \Omega^1_Y(\log S)))=0.
$$
The equality on the right hand side implies that
$$
6\cdot \ch_2(\sO_Y\oplus \Omega^1_Y(\log S)) - 2 \cdot \ch_1(\sO_Y\oplus \Omega^1_Y(\log S))^2
\Delta(\sO_Y\oplus \Omega^1_Y(\log S))=0.
$$
(see \cite[Lemma 3.2]{VZ4}, for example). So
$$
3\cdot \ch_2(\Omega^1_Y(\log S)) = \ch_1(\Omega^1_Y(\log S))^2
$$
and $U$ is a ball quotient by Yau's uniformization theorem.

For the last step, recall that
$$
\ch_1(G)= 6\cdot \ch_1(F^{1,0}) + 4\cdot \ch_1(\omega_Y(S))=0.
$$
So modulo ${\rm Pic}_0(Y)$ the invertible sheaves $\det(F^{1,0})^3$ and $\omega_Y(S)^2$ coincide, and
$\omega_Y(S)$ is $3$-divisible, hence $\omega_Y(S)^{\frac{1}{3}}$ exists.
\begin{equation}\label{eq.11}
(\omega_Y(S)^{\frac{-1}{3}}\otimes \Omega^1_Y(\log S) \oplus \omega_Y(S)^{\frac{-1}{3}} \ , \ {\rm id})
\end{equation}
is a polystable Higgs bundle of degree zero. Note that $Y$ is a smooth Mumford compactification of $U$, and 
that (\ref{eq.11}) is the universal Hodge bundle on the period domain
$B^2=\tilde U$. Hence, it is induced by the local system $\L$ corresponding to the discrete subgroup
$$
\pi_1(U)=\Gamma\subset SU(2,1).
$$
Obviously $\L$ is irreducible of rank $3$, and $\L$ is not invariant under complex conjugation. Since the
Higgs bundles coincide, one finds
\begin{equation}\label{eq.12}
\bigwedge^3(R^1f_*\C_V)=\bigwedge^3(\L\oplus \L^\vee).
\end{equation}
Let $\rho':\pi_1(Y)\to \Sp(6,\R)$ and $\rho:\pi_1(Y)\to \Sp(6,\R)$ be the representations giving
$R^1f_*\C_V$ and $\L\oplus \L^\vee$, respectively. Since their third wedge products coincide, they only differ by a constant, up to conjugation. So $R^1f_*\C_V$ is also the direct sum of two
rank $3$ local systems, say $\W_1$ and $\W_2$. The equality (\ref{eq.12}) can be rewritten as
$$
S^2(\W_1)\oplus S^2(\W_2)\oplus \C^2=
S^2(\L)\oplus S^2(\L^\vee)\oplus\C^2.
$$
The direct factors on the left hand side are all irreducible, hence
up to renumbering, one finds $S^2(\W_1)=S^2(\L)$. This implies that over an \'etale covering of $Y$
one has $\W=\L$, as claimed in i).

Assuming i) the equality of the Chern numbers in is easy to obtain. In fact, by Example \ref{su.1}
one has a unitary invertible sheaf $\sU$ on $Y$ with
$$
F^{1,0}=\omega_Y(S)^{\frac{-1}{3}}\otimes \Omega^1_Y(\log S)\otimes \sU \oplus
\omega_Y(S)^{\frac{1}{3}}\otimes \sU^\vee,
$$
hence
\begin{multline*}
\ch_1(f_*\omega_{X/Y})=\ch_1(F^{1,0})=\\
\frac{-2}{3}\cdot\ch_1(\omega_Y(S)) + \ch_1(\Omega^1_Y(\log S)) +
\frac{1}{3}\cdot\ch_1(\omega_Y(S))=\frac{2}{3}\cdot\ch_1(\omega_Y(S)).
\end{multline*}
\end{proof}
For $g>3$ it is hard to control the largest $r$ with $\tau^r\neq 0$.
Fixing $r$, one obtains:
\begin{proposition}\label{gy.3}
Assume in Proposition \ref{gy.1} that the Arakelov equality (\ref{eq.8}) holds for
all irreducible $\C$ subvariations of Hodge structures in $R^1f_*\C_V$, and
that $\tau^{r+1}=0$. Then
$$
\ch_1(f_*\omega_{X/Y}). \ch_1(\omega_Y(S)) \leq \frac{r}{3}\ch_1(\omega_Y(S))^2.
$$
\end{proposition}
\begin{proof}
By \cite[Theorem 1]{VZ4} the Arakelov equality implies that $F^{1,0}$ and $F^{0,1}$ are both $\mu$-semistable, and that
$$
\tau:F^{1,0}\to F^{0,1}\otimes \Omega^1_Y(\log S)
$$
is a morphism between $\mu$-semistable sheaves of the same slope. Writing
$$
T_Y(-\log S)= \Omega^1_Y(\log S)^\vee,
$$
the same holds  for the induced map
$$
\tau'_q :\det(F^{1,0})\otimes S^q(T_Y(-\log S)) \>>> \bigwedge^{g-q}F^{1,0}\oplus \bigwedge^q F^{0,1}.
$$
The left hand side is stable, and $G^{g-q,q}$, as the image of $\tau'_q$, is either zero or isomorphic to
$$
\det(F^{1,0})\otimes S^q(T_Y(-\log S)).
$$
Choosing $r$ to be minimal in \ref{gy.3} we may assume that it is the largest integer with $G^{g-r,r}\neq 0$. Since $G$ is a Higgs subbundle of the Higgs bundle of a local constant system, by \ref{ar.3}
$$
\ch_1(G).\ch_1(\omega_Y(S))\leq 0.
$$
The sheaf $G^{g-q,q}$ has slope $\mu(\det(F^{1,0}))-q\cdot\mu(\Omega^1_Y(\log S))$ and
\begin{multline*}
\ch_1(G).\ch_1(\omega_Y(S))=\sum_{q=0}^r (q+1)\cdot\big(\mu(\det(F^{1,0})) - q\cdot\mu(\Omega^1_Y(\log S))\big)=\\
\frac{(r+1)\cdot(r+2)}{2}\cdot \mu(\det(F^{1,0})) -\\
\big(\frac{r\cdot(r+1)\cdot(2r+1)}{6}+\frac{r\cdot(r+1)}{2}\big) \cdot \mu(\Omega^1_Y(\log S)).
\end{multline*}
The latter is a positive multiple of
$$
\mu(\det(F^{1,0}))-\frac{2\cdot r}{3}\cdot \mu(\Omega^1_Y(\log S)).
$$
\end{proof}
For $r=g-1$ and $g\geq 4$ the inequality in Proposition \ref{gy.3} is worse than the Arakelov inequality (\ref{eq.10}). In \cite[Section 9]{VZ4} we calculated the maximal $r$ with
$\tau^r\neq 0$ for the ball quotients. In particular for the Picard modular surfaces
and for the variation of Hodge structures in Theorem \ref{su.2}, b), one obtains
$r=2\cdot \rk(\U)=\frac{2\cdot g}{3}$. Then Proposition \ref{gy.3} implies that
\begin{equation}\label{eq.13}
\ch_1(f_*\omega_{X/Y}). \ch_1(\omega_Y(S)) \leq \frac{2\cdot g}{9}\ch_1(\omega_Y(S))^2.
\end{equation}
We leave it to the reader to verify, that the Arakelov equality for $\L\otimes\U$ and for its dual 
are just saying that (\ref{eq.13}) is an equality.

\end{document}